\documentclass[11pt]{article}

\usepackage{times,amsfonts}
\usepackage[a4paper]{geometry}
\newtheorem{thm}{Theorem}[section]
\newtheorem{lem}[thm]{Lemma}
\newenvironment{proof}{\noindent {\bf Proof}.}{\hfill$\Box$\\[-5mm]}

\newcommand{\NN}{\mathbb{N}}
\newcommand{\RR}{\mathbb{R}}
\newcommand{\DD}{\mathbb{D}}

\newcommand{\EE}{\mathbf{E}}
\newcommand{\PP}{\mathbf{P}}

\newcommand{\One}{{\mathrm{1} \kern -0.27em \mathrm{I}}}

\newcommand{\defeg}{\mathrel :=}

\newcommand{\ut}{{\tilde u}}
\newcommand{\lin}{\mathbf{K}}
\newcommand{\lint}{\mathcal{K}}

\begin{document}

\title{A functional central limit theorem in equilibrium 
for a large network in which customers join the shortest of several queues}

\author{\sc Carl Graham~\footnote{CMAP, 
{\'E}cole Polytechnique, 91128 Palaiseau, France. UMR CNRS 7641. 
\tt carl@cmapx.polytechnique.fr}
}
\date{}

\maketitle

\begin{small}
\noindent
September 9, 2003.

\noindent
{\bf Abstract}.
We consider
$N$ single server infinite buffer queues with service rate $\beta$.
Customers arrive at rate $N\alpha$,
choose $L$ queues uniformly, and join the shortest one.
The stability condition is $\alpha < \beta$.
We study in equilibrium
the fraction of queues of length at least $k\ge0$.
We prove a functional central limit theorem 
on an infinite-dimensional Hilbert space with its weak topology,
with limit a stationary Ornstein-Uhlenbeck process.
We use ergodicity  and
justify the  inversion of limits 
$\lim_{N\to\infty}\lim_{t \to\infty} =\lim_{t\to\infty}\lim_{N \to\infty}$
by a compactness-uniqueness method.
The main tool for proving tightness of the ill-known invariant laws
and ergodicity of the limit
is a global  exponential stability result
for the nonlinear dynamical system obtained in the functional law of large numbers
limit.


\end{small}
\bigskip\hrule

\section{Introduction}

\subsection{The queuing network, and some notation}

Customers arrive at rate $N\alpha$ on a network 
constituted of $N\ge L\ge1$ infinite buffer single server 
queues.
Each customer is allocated
$L$ distinct queues uniformly at random and joins the shortest, ties
being resolved uniformly. 
Servers work at rate $\beta$.
Inter-arrival times, allocations, and services 
are independent and memoryless.
For $L=1$ we have $N$ i.i.d.
$M_\alpha/M_\beta/1/\infty$ queues, and for $L\ge2$
the interaction structure depends only on sampling 
from the empirical measure of $L$-tuples of queue
states. In statistical mechanics terminology, this system is in 
$L$-body mean-field interaction.

The process $(X^N_i)_{1\leq i\leq N }$, where
$X^N_i(t)$ denotes the length of queue $i$ at time $t\ge0$, is Markov.
Its empirical measure $\mu^N$ 
with samples in $\mathcal{P}(\DD(\RR_+,\NN) )$
and its marginal process 
$\bar X^N = (\bar X^N_t)_{t\ge0}$ 
with sample paths
in $\DD(\RR_+,\mathcal{P}(\NN))$
are given by
\[
\mu^N={1 \over N}\sum_{i=1}^N \delta_{X^N_i}\,,
\qquad
\bar X^N_t
= {1 \over N}\sum_{i=1}^N \delta_{X^N_i(t)}
\,.
\] 

We are interested in the tails of the distributions $\bar X^N_t$. 
We consider
\[
\mathcal{V} 
= \Bigl\{ 
(v(k))_{k\in\NN} : v(0)=1,\ v(k)\geq v(k+1),\ \lim_{k\to\infty} v(k) =0 
\Bigr\} 
\subset c_0\,,
\quad 
\mathcal{V}^N = \mathcal{V} \cap {1 \over N} \NN^\NN\,,
\]
with the uniform topology. Note that the uniform and the product topology
coincide on $\mathcal{V}$. We consider the 
process $R^N=(R^N_t)_{t\ge0}$
with sample paths in $\DD(\RR_+,\mathcal{V}^N)$
 given by 
\[
R^N_t(k) = {1 \over N}\sum_{i=1}^N \One_{X^N_i(t)\ge k}\,,
\qquad k\in\NN\,,
\] 
the fraction of queues at time $t$ of length 
at least $k$.

We have $R^N_t(k) = \bar X^N_t([k,\infty [)$ 
and $\bar X^N_t\{k\} = R^N_t(k) - R^N_t(k+1)$ using the classical  homeomorphism
between $\mathcal{P}(\NN)$ and  $\mathcal{V}$, which
maps the subspace of probability measures with 
finite first moment
onto  $\mathcal{V}\cap \ell_1$ corresponding to having 
a finite number of customers.
The symmetry structure implies that $\bar X^N$ and $R^N$
are Markov processes.

The network is ergodic if and only if $\alpha < \beta$
(Theorem~5~(a) in \cite{Vved:96},
Theorem~4.2 in \cite{Graham:00}). The proofs use non-constructive 
ergodicity criteria, and
we lack information and controls on the
invariant laws (stationary distributions). 
We study the large $N$ asymptotics in the stationary regime using an
indirect approach involving ergodicity in appropriate transient regimes
and an inversion of limits 
for large $N$ and large times.
Law of large numbers (LLN) results are already known, and
we shall obtain a functional central limit theorem (CLT).

\noindent
{\em General notation\/}.
We denote by $c^0_0$ and $\ell^0_p$ for $p\ge1$
the subspaces  of sequences vanishing at $0$
of the classical sequence spaces $c_0$ (with limit 0) and
$\ell_p$ (with summable $p$-th power).
The diagonal matrix with successive diagonal terms given by
the sequence $a$ is denoted by $\mathrm{diag}(a)$.
When using matrix notations, sequences vanishing at $0$
are often identified with infinite column vectors
indexed by $\{1,2,\cdots\}$.
Sequence inequalities, etc., should be interpreted termwise.
Empty sums are equal to
$0$ and empty products to $1$. Constants such as $K$ may
vary from line to line.
We denote by $g_\theta = (\theta^k)_{k\ge1}$ the geometric
sequence of reason $\theta$.

\subsection{Previous results: laws of large numbers}

We relate results found in essence in
Vvedenskaya et al.~\cite{Vved:96}.
Graham~\cite{Graham:00} extended some of these results, and also
considered  the  empirical measures on path space 
$\mu^N$, yielding chaoticity results (asymptotic independence of queues).
(The rates $\nu$ and $\lambda$
in \cite{Graham:00} are replaced here by $\alpha$ and $\beta$.)

Consider the  mappings with values in $c^0_0$ given
for $v$ in $c_0$ by
\begin{equation}
\label{F}
F_+(v)(k) = \alpha\!\left( v(k-1)^L -v(k)^L  \right),
\quad 
F_-(v) (k) = \beta  (v(k) - v(k+1))\,,
\qquad 
k\ge1\,,
\end{equation}
and $F = F_+ - F_-$ and the nonlinear differential equation 
$\dot u = F(u)$
on $\mathcal{V}$ given for $t\ge0$ by 
\begin{equation}
\label{is}
\dot u_t(k)
= F(u_t)(k)
= \alpha\! \left( u_t(k-1)^L-u_t(k)^L \right)
- \beta\! \left( u_t(k) - u_t(k+1)\right)\,,
\qquad k\geq 1\,.
\end{equation}
This is the infinite system of scalar
differential equations
(1.6) in \cite{Vved:96} (where the arrival rate is $\lambda$ and
service rate $1$) and  (3.9) in  \cite{Graham:00}.
Note that $F_-$ is linear.

\begin{thm}
\label{euv}
There  exists a unique solution $u=(u_t)_{t\ge0}$
taking values in $\mathcal{V}$ for (\ref{is}), and $u$ is in $C(\RR_+,\mathcal{V})$. 
If $u_0$ is in $\mathcal{V}\cap\ell_1$
then $u$ takes values in $\mathcal{V}\cap\ell_1$.
\end{thm}

\begin{proof}
We use Theorem~3.3  and Proposition 2.3 in \cite{Graham:00}. 
These exploit the homeomorphism
between $\mathcal{P}(\NN)$ with the weak topology 
and  $\mathcal{V}$ with the product topology.
Then (\ref{is}) corresponds to a non-linear forward Kolmogorov equation for a pure jump
process with uniformly bounded (time-dependent) jump rates.
Uniqueness within the class of bounded measures
and existence of a probability-measure valued solution 
are obtained using the total variation norm.
Theorem~1~(a) in \cite{Vved:96} yields existence (and uniqueness) 
in $\mathcal{V}\cap\ell_1$.
\end{proof}

Firstly, a functional LLN
for initial conditions satisfying a LLN
is part of Theorem~3.4 in \cite{Graham:00} and can be deduced from
Theorem~2 in \cite{Vved:96}.

\begin{thm}
\label{lln}
Assume that
$(R^N_0)_{N\ge L}$ converges in law  to
$u_0$ in $\mathcal{V}$.
Then  $(R^N)_{N\ge L}$ converges in law in $\DD(\RR_+,\mathcal{V})$
to the unique solution $u=(u_t)_{t\ge0}$ 
starting at $u_0$ for (\ref{is}).
\end{thm}

Secondly, the limit equation  (\ref{is}) has a globally attractive stable point $\ut$ 
in $\mathcal{V}\cap \ell_1$. 

\begin{thm}
\label{utilde}
For $\rho = \alpha / \beta <1$ the  equation (\ref{is}) has
a unique stable point $\ut$ in $\mathcal{V}$ given by
\[
\ut = (\ut(k)_{k\in\NN}\,,
\qquad
\ut(k)= \rho^{(L^k -1) / (L-1)}= \rho^{L^{k-1}+L^{k-2}+\cdots + 1}\,,
\]
and the solution $u$ of (\ref{is}) 
starting at any $u_0$ in $\mathcal{V}\cap \ell_1$  is such that
$\lim_{t\to\infty} u_t = \ut$.
\end{thm}
\begin{proof}
Theorem~1~(b) in \cite{Vved:96} yields that $\ut$ 
is  globally asymptotically stable in $\mathcal{V}\cap \ell_1$.
A stable point $u$ in $\mathcal{V}$ satisfies
$
\beta u(k+1) - \alpha u(k)^L = \beta u(k) - \alpha u(k-1)^L
= \cdots =
\beta u(1) - \alpha 
$
and converges to $0$, hence
$u(1) = \alpha/\beta$ and
$u(2)$,  $u(3)$, \dots\, are successively determined uniquely.
\end{proof}

Lastly, a compactness-uniqueness method 
justifying the inversion of limits 
$\lim_{N\to\infty}\lim_{t\to\infty} = \lim_{t\to\infty}\lim_{N\to\infty}$
yields a result in equilibrium. This 
method was
used by Whitt~\cite{Whitt:85} for the star-shaped loss network, and
is described in detail in Graham~\cite{Graham:99} Sections~9.5 and 9.7.3.
The following functional LLN in equilibrium (Theorem~4.4 in \cite{Graham:00})
can be deduced from \cite{Vved:96}, but is not stated there as such; it implies
using uniform integrability bounds that under the invariant laws
$\lim_{N\to\infty} \EE( R^N_0(k) ) = \ut(k)$ for $k\in\NN$, a result stated in 
Theorem~5 (c) in \cite{Vved:96}.

\begin{thm}
\label{lln.eq}
Let $\rho = \alpha / \beta <1$ and
the networks of size $N\ge L$ be in equilibrium. Then $(R^N)_{N\ge L}$
converges in probability in $\DD(\RR_+,\mathcal{V})$ to $\ut$.
\end{thm}

Note that $\ut(k)$ decays hyper-exponentially in $k$ for $L\ge2$ 
instead of the exponential decay $\rho^{k}$ corresponding to 
i.i.d. queues in equilibrium ($L=1$).
The asymptotic large queue sizes are dramatically decreased 
by this simple choice.

We seek rates of convergence and confidence intervals. 
Theorem~3.5 in  \cite{Graham:00} gives convergence bounds 
when $(X^N_i(0))_{1\leq i \leq N}$ are i.i.d.\ for 
the variation norm 
on ${\cal P}(\DD([0,T],\NN^k))$ 
using results in Graham and M{\'e}l{\'e}ard~\cite{Meleard:94}. This can be extended 
if the initial laws satisfy a priori controls, 
but it is not so in equilibrium, where on the contrary controls are obtained 
using the network evolution.

\subsection{The outline of this paper}

We consider the process $R^N$ with values in  $\mathcal{V}^N$, a
solution $u=(u_t)_{t\ge0}$ for (\ref{is}) in $\mathcal{V}$, and
the empirical fluctuation processes $Z^N = (Z^N_t)_{t\ge0}$
with sample paths in $c^0_0$ given by
\begin{equation}
\label{flu}
 Z^N = N^{1/2} (R^N - u)\,,
\qquad
Z^N_t = N^{1/2} (R^N_t - u_t)\,.
\end{equation}
We are interested in particular in the stationary regime, which defines {\em implicitly\/}
the initial data:
the law of $R^N_0$ is the invariant law for $R^N$
and $u_0=\ut$. 

Our main result is a functional CLT: in equilibrium
$(Z^N)_{N\ge L}$ converges in law to 
a stationary Ornstein-Uhlenbeck process, which we characterize.
This {\em implies\/} a CLT for the marginal laws:
under the invariant laws $(Z^N_0)_{N\ge L}$ converges
to the invariant law for this Gaussian process. This
important result seems very difficult to obtain 
directly. We use ergodicity of $Z^N$ for fixed $N$ 
and intricate fine studies of the long-time behavior of the  
nonlinear dynamics  appearing at the 
large $N$ limit, simply in order to prove
tightness bounds for  $(Z^N_0)_{N\ge L}$ under the invariant laws and 
ergodicity for the Ornstein-Uhlenbeck process. 

Section~\ref{clt} introduces the main theorems, which are proved in subsequent sections.
Section~\ref{derOU} considers arbitrary $u_0$ and $R^N_0$ and derives
martingales of interest and the limit Ornstein-Uhlenbeck process. 
We consider the stationary regime whenever possible for simplicity,
but the infinite-horizon bounds used for the control of the 
invariant laws are obtained considering {\em transient\/} regimes.

We study the Ornstein-Uhlenbeck process in Section~\ref{propOU}.
We give a spectral representation for the linear operator in the drift term,
and prove the existence of a spectral gap. A main difficulty is that
the Hilbert space in which this operator is self-adjoint
is {\em not\/} large enough (its norm is {\em too\/} strong)
for the limit non-linear dynamical system 
and for the invariant laws for finite $N$.
We obtain results of global exponential stability
in appropriate Hilbert spaces in which it is {\em not\/} self-adjoint.

In Section~\ref{expstab} we  prove that $\ut$ is globally
exponentially  stable for the non-linear 
dynamical system in appropriate Hilbert spaces.
In Section~\ref{cfp},  uniformly for large  $N$, we obtain bounds 
for the processes $Z^N$ on $[0,T]$ 
using martingale properties, and then for  $Z^N_t$ uniformly for $t\ge0$
using  the above result on the dynamical system in order 
to  iterate the bounds on intervals of length $T$. 
Bounds on the invariant laws of $Z^N$ follow using ergodicity.
We then prove
the functional CLT by a 
compactness-uniqueness method
and martingale characterizations. We 
consider the non-metrizable weak topology on the Hilbert spaces, 
and use adapted tightness criteria
and the above bounds.

\section{The functional central limit theorem in equilibrium}
\setcounter{equation}{0}
\label{clt}

In this paper we concentrate on the stationary regime, and
assume that $\rho = \alpha/\beta <1$ and $u_0=\ut = u$.
We leave the explicit study of transient regimes for a forthcoming paper.
We quickly introduce notation and state the main results,
leaving most proofs for later.

\subsection{Preliminaries}

For any sequence $w=(w(k))_{k\ge 1}$ such that $w>0$
we define the Hilbert spaces
\[
L_2(w) = \biggl\{ x \in \RR^\NN : x(0)=0\,,\ 
\Vert x \Vert_{L_2(w)}^2 = \sum_{k\ge1} \biggl({x(k)\over w(k)}\biggr)^2 w(k)
= \sum_{k\ge1} x(k)^2 w(k)^{-1} <\infty
\biggr\}
\]
and in matrix notation $(x,y)_{L_2(w)} = x^* \mathrm{diag}(w^{-1})y$.
We  consider the elements of
$L_2(w)$ as measures identified with 
their densities with respect to the reference measure $w$. 
Then $L_1(w) = \ell^0_1$ and if
$w$ is summable then 
$\Vert x \Vert_1 \le \Vert w \Vert_1^{1/2}\Vert x \Vert_{L_2(w)}$
and $L_2(w) \subset \ell^0_1$.
Using $L_2(1) = \ell^0_2$ as a pivot space,
for bounded $w$ we have the Gelfand triplet of Hilbert spaces 
$
L_2(w) \subset \ell^0_2  \subset L_2(w)^* = L_2(w^{-1})
$.

\begin{lem}
\label{eqiv.l2w}
If $w = O(v)$ and $v = O(w)$
then the $L_2(v)$ and  $L_2(w)$ norms
are equivalent.
\end{lem}
\begin{proof}
This follows from obvious computations.
\end{proof}

We give a refined existence result for $(\ref{is})$.  We recall that
$g_\theta = (\theta^k)_{k\ge1}$.

\begin{thm}
\label{isw}
Let $w>0$ be such that there exists  $c>0$ and $d>0$ with
\[
c  w(k+1)\le  w(k) \le d w(k+1)\,,
\qquad k\ge1\,.
\]
Then in $\mathcal{V} \cap L_2(w)$ the mappings
$F$, $F_+$ and $F_-$ are Lipschitz for 
the $L_2(w)$ norm and there is existence and uniqueness
for (\ref{is}).
The assumptions and conclusions hold for $w = g_\theta$ for $\theta>0$.
\end{thm}
\begin{proof}
The identity $x^L-y^L = (x-y)(x^{L-1} + x^{L-2}y +\cdots + y^{L-1})$ yields
\begin{eqnarray*}
\left(u(k-1)^L - v(k-1)^L\right)^2 w(k)^{-1}
&\le& 
\left(u(k-1) - v(k-1)\right)^2 L^2 d w(k-1)^{-1}\,,
\\
\left(u(k)^L - v(k)^L\right)^2 w(k)^{-1}
&\le& 
\left(u(k) - v(k)\right)^2 L^2  w(k)^{-1}\,,
\\
\left(u(k+1) - v(k+1)\right)^2  w(k)^{-1}
&\le&
\left(u(k+1) - v(k+1)\right)^2  c^{-1}w(k+1)^{-1}\,,
\end{eqnarray*}
hence we have the Lipschitz bounds
$
\Vert F_+(u) - F_+(v) \Vert_{L_2(w)}^2 
\le 2\alpha^2 L^2 (d+1)
\Vert u - v  \Vert_{L_2(w)}^2
$ and
$
\Vert F_-(u) - F_-(v) \Vert_{L_2(w)}^2 
\le 2\beta^2  (c^{-1}+1)
\Vert u - v  \Vert_{L_2(w)}^2
$
and existence and uniqueness follows by a classical
Cauchy-Lipschitz method.
We have
$ \theta^{-1} \theta^{k+1} \le \theta^k \le \theta^{-1} \theta^{k+1}$
for $k\ge1$.
\end{proof}

\subsection{The Ornstein-Uhlenbeck process}
\label{oueq}

We consider the linear operator 
$\lint : x \in c_0^0 \mapsto \lint x \in c_0^0$
given  by
\begin{eqnarray}
\lint x(k) &=&
\alpha L {\ut}(k-1)^{L-1} x(k-1)
-\left(\alpha L {\ut}(k)^{L-1} + \beta\right) x(k)
+  \beta x(k+1)
\nonumber\\
\label{defDst}
&=&
\beta L \rho^{L^{k-1}} x(k-1)
-\left(\beta L \rho^{L^k} + \beta\right) x(k)
+  \beta x(k+1)\,,
\qquad k\ge1\,,
\end{eqnarray}
which we identify
with its infinite matrix in the canonical basis
$(0,1,0,0 \ldots), (0,0,1,0 \ldots), \dots$ 
\begin{equation}
\label{matrixK}
\lint =
\pmatrix{
 -\left(\beta L \rho^{L} + \beta\right) &  \beta & 0  & 0 &\cdots 
\cr
 \beta L \rho^{L} & -\left( \beta L \rho^{L^2} + \beta\right) &  \beta   & 0 & \cdots
\cr
 0 &  \beta L \rho^{L^2} & -\left( \beta L \rho^{L^3} + \beta\right)  & \beta & \cdots
\cr
 0 &  0 &  \beta L \rho^{L^3} & -\left( \beta L \rho^{L^4} + \beta\right)  &  \cdots
\cr
 \vdots & \vdots & \vdots  & \vdots & 
}
\end{equation}
used identifying the sequence $x=(0,x(1),x(2), \ldots\,)$  with
its coordinates  in the canonical basis $(x(1),x(2), \ldots\,)$ 
taken as a column vector.

Note that $\lint=\mathcal{A}^*$ where $\mathcal{A}$  is  
the infinitesimal generator of a sub-Markovian birth and death process.
We shall develop this point of view and obtain a spectral decomposition for $\lint$ 
in Section~\ref{spectdec}, to which we give a few anticipated 
references below. 
The potential coefficients  of $\mathcal{A}$ given by 
\[
\pi = (\pi(k))_{k\ge1}\,,
\qquad
\pi(k) =   L^{k-1} \rho^{(L^k-L)/(L-1)} =  \rho^{-1} L^{k-1} {\ut}(k)\,,
\]
solve the detailed balance equations $\pi(k+1) = L\rho^{L^k}\pi(k)$
with $\pi(1)=1$. 

The linearization of (\ref{is}) around its stable point $\ut$
is the linearization of the equation satisfied by $z=u-\ut$ and is 
given for $t\ge0$ by the forward Kolmogorov equation
\begin{equation}
\label{linis}
\dot z_t = \lint z_t\,.
\end{equation} 

Let $B=(B(k))_{k\in\NN}$ be independent Brownian motions such that
$B(0)=0$ and $\mathrm{var}(B_1(k)) = \EE(B_1(k)^2) = \tilde v(k)$ 
where $\tilde v$ in $c^0_0$ is given by
\[
\tilde v(k)
= 2 \beta \left(\ut(k) - \ut(k+1)\right)
= 2 \beta  \rho^{(L^k-1)/ (L-1)} \left(1 - \rho^{L^{\smash{k}}} \right),
\qquad k\ge1\,.
\]
The  infinitesimal covariance matrix of $B$ is given by
$\mathrm{diag}(\tilde v)$. 

\begin{thm}
The process $B$ is an Hilbertian Brownian motion in $L_2(w)$ 
if and only if
\begin{equation}
\label{chbm}
\sum_{k\ge1} \ut(k) w(k)^{-1}
=
\sum_{k\ge1} \rho^{(L^k-1)/ (L-1)} w(k)^{-1} <\infty\,.
\end{equation}
This is true for 
$w= \pi$ and $w = g_\theta$ for $\theta>0$ when $L\ge2$ or 
for $w = g_\theta$ for $\theta>\rho$ when $L=1$.
\end{thm}
\begin{proof}
This follows from obvious computations.
\end{proof}
\goodbreak

The Ornstein-Uhlenbeck process $Z = (Z(k))_{k\in\NN}$ solves
the affine SDE given for $t\ge0$ by 
\begin{equation}
\label{sde}
Z_t = Z_0 + \int_0^t \lint Z_s \,ds + B_t
\end{equation} 
which is a Brownian perturbation of (\ref{linis}).

\begin{thm}
\label{Dbdd}
 Let $w>0$ be such that there exists  $c>0$ and $d>0$ with
\[
c  w(k+1)\le  w(k) \le d \rho^{-2L^k} w(k+1)\,,
\qquad k\ge1\,.
\]
(a)
In $L_2(w)$, the operator $\lint$ is bounded,
equation (\ref{linis}) has a unique solution 
$z_t = \mathrm{e}^{\lint t} z_0$ where
$\mathrm{e}^{\lint t}$  has a spectral representation 
given by (\ref{repr}), and
there is uniqueness of solutions for the SDE (\ref{sde}).
The assumptions and  conclusions hold 
for $w= \pi$ and $w = g_\theta$ for $\theta>0$.
\\
(b) In addition let $w$ satisfy (\ref{chbm}).
The SDE (\ref{sde}) has a unique solution
$
Z_t = \mathrm{e}^{\lint t} Z_0 
+ \int_0^t \mathrm{e}^{\lint(t-s)}\,dB_s
$ in $L_2(w)$, further explicited in (\ref{explou}).
The assumptions and conclusions hold  for 
$w= \pi$ and $w = g_\theta$ for $\theta>0$ when $L\ge2$ or
for $w = g_\theta$ for $\theta>\rho$ when $L=1$ .
\end{thm}

\begin{thm}
\label{spe.gap}
(Spectral gap.)
The operator $\lint$ is bounded self-adjoint in $L_2(\pi)$.
The least point $\gamma$ of the spectrum of $\lint$
is such that $0  < \gamma \le \beta$. The 
solution $z_t = \mathrm{e}^{\lint t}z_0$ for (\ref{linis}) in $L_2(\pi)$
satisfies
$\Vert z_t \Vert_{L_2(\pi)} 
\le \mathrm{e}^{-\gamma t}  \Vert z_0 \Vert_{L_2(\pi)}$.
\end{thm}

The  $L_2(\pi)$ norm is too strong for studying the CLT. Indeed, 
$
\PP(X^N_1 + \cdots + X^N_N \ge N k) \le 
\PP(X^N_1\ge k) + \cdots + \PP(X^N_N \ge k)
$ 
and since the total service rate in the system cannot exceed $N\beta$,
by comparison
with an $M_{N\alpha}/M_{N\beta}/1$ queue, in equilibrium
\[
\EE(R^N_t(k)) = \PP(X^N_i(t) \ge k) \ge {1\over N} \rho^{Nk}
\]
decreases at most
exponentially in $k\ge0$.
Further,
the mapping $F_+$ is not Lipschitz in 
$\mathcal{V} \cap L_2(\pi)$ for 
the $L_2(\pi)$ norm,
see Theorem~\ref{isw} and
the contrasting assumptions and proof of Theorem~\ref{Dbdd}.
We prove global exponential stability in appropriate spaces.  

\begin{thm}
\label{linis.exp}
Let  $0 < \theta < 1$ when $L\ge2$ or $\rho \le \theta < 1$ when $L=1$. 
There  exists  $\gamma_\theta>0$ and $C_\theta < \infty$ such that
the solution $z_t = \mathrm{e}^{\lint t}z_0$ for (\ref{linis})  
in $L_2(g_\theta)$
satisfies
$\Vert z_t \Vert_{L_2(g_\theta)} 
\le \mathrm{e}^{-\gamma_\theta t} C_\theta \Vert z_0 \Vert_{L_2(g_\theta)}$.
\end{thm}

We deduce exponential ergodicity
for the Ornstein-Uhlenbeck process,
valid for any $w$ satisfying the conclusions of 
Theorems \ref{Dbdd} and \ref{linis.exp}.

\begin{thm}
\label{ou.erg}
Let $w= \pi$ or $w = g_\theta$ with $0 < \theta<1$ when $L\ge2$ or
let $w = g_\theta$ with $\rho < \theta<1$ when $L=1$.
Any solution for the SDE (\ref{sde}) in $L_2(w)$
converges in law for large times 
 to its unique invariant law (exponentially fast). This law is
the law of $\int_0^\infty \mathrm{e}^{\lint t}dB_t$ which is 
 Gaussian centered with covariance matrix 
$\int_0^\infty \mathrm{e}^{\lint t} \mathrm{diag}(\tilde v) \mathrm{e}^{\lint^* t}dt$,
further explicited in  (\ref{limlaw}) and (\ref{cov}).
There is a unique stationary Ornstein-Uhlenbeck process
solving the SDE (\ref{sde}) in  $L_2(w)$.
\end{thm}

\subsection{Global exponential stability for the dynamical system 
and tightness estimates}

Global  exponential stability of the dynamical system allows 
control of the invariant laws
using the long time behavior.
We need  uniformity over the state space, and
Theorems \ref{spe.gap} or \ref{linis.exp} are useless for this purpose
(except in the linear case $L=1$).
Such a result does {\em not\/} hold in $L_2(\pi)$ for $L\ge2$.

\begin{thm}
\label{glexst}
Let $\rho \le \theta < 1$
and
$u$ be the solution of (\ref{is}) starting at 
$u_0$ in $\mathcal{V} \cap L_2(g_\theta)$.
There exists  $\gamma_\theta>0$ and $C_\theta < \infty$
such that
$\Vert u_t -\ut \Vert_{L_2(g_\theta)} 
\le  \mathrm{e}^{- \gamma_\theta t} C_\theta \Vert u_0 -\ut \Vert_{L_2(g_\theta)}$.
\end{thm}

The following finite-horizon bounds yield
tightness estimates for the processes $(Z^N)_{N\ge L}$
provided the initial laws are known to satisfy similar bounds.

\begin{lem}
\label{l2bds}
For $\theta >0$ and $T\ge0$ we have
\[
\limsup_{N\ge L} 
\EE\left(\left\Vert Z^N_0 \right\Vert_{L_2(g_\theta)}^2\right) < \infty
\Rightarrow
\limsup_{N\ge L} 
\EE\biggl( \sup_{0 \le t\le T} \left\Vert Z^N_t \right\Vert_{L_2(g_\theta)}^2\biggr) 
< \infty\,.
\]
\end{lem}
Theorem~\ref{glexst} is an essential ingredient in the proof of
the following infinite-horizon bound
for the marginal laws of the processes. 

\begin{lem}
\label{l2bdsinf}
Let  $\rho \le \theta < 1$ when $L\ge2$ or  $\rho < \theta < 1$ when $L=1$. Then
\[
\limsup_{N\ge L} 
\EE\left(\left\Vert Z^N_0 \right\Vert_{L_2(g_\theta)}^2\right) < \infty
\Rightarrow
\limsup_{N\ge L} \sup_{t\ge0}
\EE\left(\left\Vert Z^N_t \right\Vert_{L_2(g_\theta)}^2\right) < \infty\,.
\]
\end{lem}
This yields control of the long time limit of the marginals,
the  invariant law, which in turn will enable us to use Lemma~\ref{l2bds}
to prove tightness of the processes in equilibrium.
\begin{lem}
\label{tight.in}
Let  $\rho \le \theta < 1$ when $L\ge2$ or  $\rho < \theta < 1$ when $L=1$. Then
under the invariant laws
\[
\limsup_{N\ge L} \EE\left(\Vert Z^N_0 \Vert_{L_2(g_\theta)}^2\right) < \infty\,.
\]
\end{lem}

\subsection{The main result: the functional CLT in equilibrium}

This result is obtained by a 
compactness-uniqueness method.
We refer to Jakubowski~\cite{Jakubowski:86} for the Skorokhod topology
for the non-metrizable weak topology on 
infinite-dimensional Hilbert spaces.

\begin{thm}
\label{eqfclt}
Let the networks of size $N\ge L$ be in equilibrium.
For $L\ge2$ consider 
$L_2(g_{\rho})$
with its weak topology and
$\DD(\RR_+, L_2(g_{\rho}))$ with the corresponding Skorokhod
topology.
 Then $(Z^N)_{N\ge L}$
converges in law to the unique stationary Ornstein-Uhlenbeck process 
solving the SDE~(\ref{sde}),
which is  continuous
and Gaussian, in particular  $(Z^N_0)_{N\ge L}$ converges in law
to the invariant law for this process (see  Theorem~\ref{ou.erg}).
For $L=1$ the same result holds in  $L_2(g_{\theta})$ for $\rho<\theta<1$.
\end{thm}

\section{The derivation of the limit Ornstein-Uhlenbeck process}
\setcounter{equation}{0}
\label{derOU}

Let $(x)_k = x(x-1)\cdots(x-k+1)$ for $x\in \RR$ 
denote the Jordan or falling factorial of degree $k\in \NN$.
Considering (\ref{F}),
let  the mappings $F^N$ and $F^N_+$ with values in $c^0_0$
be given for $v$ in $c_0$ by 
\[
F^N(v) = F^N_+(v) - F^{\vphantom{N}}_-(v)\,,\qquad
F^N_+(v) (k) = \alpha\, {(Nv(k-1))_L - (Nv(k))_L\over (N)_L}\,,
\quad
k\ge1\,.
\]
The process $R^N$ is  Markov on $\mathcal{V}^N$, and
when in state $r$ has jumps in its $k$-th coordinate, 
$k\ge1$, of size $1/N$  at rate $N F^N_+(r) (k)$
and size $-1/N$ at rate $N F_-(r) (k)$. 

\begin{lem}
\label{dyn}
Let 
$R^N_0$ be in  $\mathcal{V}^N$,
$u$ solve (\ref{is})
starting at $u_0$ in $\mathcal{V}$,
and $Z^N$ be given by (\ref{flu}). Then
\begin{equation}
\label{zeq}
Z^N_t  = Z^N_0
+ \int_0^t N^{1/2}\left( F^N(R^N_s) - F(u_s) \right) ds  + M^N_t
\end{equation}
defines an independent family of square integrable martingales
$M^N = (M^N(k))_{k\in\NN}$ 
independent of $R^N_0$ with Doob-Meyer brackets given by
\begin{equation}
\label{bracket}
\left\langle M^N(k) \right\rangle_t 
= \int_0^t \left\{ F^N_+(R^N_s)(k) + F^{\vphantom{N}}_-(R^N_s)(k) \right\} ds\,.
\end{equation}
\end{lem}

\begin{proof}
This follows from a classical application of the 
Dynkin formula.
\end{proof}

The first 
following combinatorial identity shows that
it is indifferent to choose the
$L$ queues with or without replacement at this level of precision.
The second one is a linearization formula.

\begin{lem}
\label{fact}
For $N\ge L$ and $a$ in $\RR$ we have
\[
A^N(a) 
\defeg
{(Na)_L \over (N)_L} - a^L 
= \sum_{j=1}^{L-1} (a-1)^j a^{L-j}  
\sum_{1\le i_1<\cdots<i_j\le L-1}
{i_1\cdots i_j \over (N-i_1)\cdots (N-i_j)}
\]
and $A^N(a) =  N^{-1} O(a)$ uniformly for $a$ in $[0,1]$.
We have $A^N(a) \le 0$ for $a$ in $\{0,N^{-1},2N^{-1},\ldots, 1\}$.
\end{lem}

\begin{proof}
We have
\[
{(Na)_L \over (N)_L} = \prod_{i=0}^{L-1}  { Na-i \over N-i}
= \prod_{i=0}^{L-1} \left(a + (a-1){ i \over N-i}\right)
\]
and by developing the product we obtain the first identity.
Direct inspection of the right-hand side of the identity   
shows that $A^N(a) = N^{-1} O(a)$ uniformly for $a$ in $[0,1]$.
For  $a$ in $\{0,N^{-1},2N^{-1},\ldots, 1\}$ the product
either is composed of terms which are positive and do not exceed $a$ or 
contains a term
equal to $0$, and hence does not exceed $a^L$.
\end{proof}

\begin{lem}
\label{lin}
For $N\ge L$ and $a$ and $h$ in $\RR$ we have
\[
B(a,h) 
\defeg
(a+h)^L - a^L - La^{L-1}h  
=\sum_{i=2}^{L} {L \choose i} a^{L-i} h^{i}
\]
with $B(a,h) =0$ for $L=1$ and $B(a,h) =h^2$
for $L=2$. For $L\ge2$ we have
$0\le B(a,h) \le h^{L} + \left(2^L -L -2\right) a h^2$
for $a$ and $a+h$ in $[0,1]$. 
\end{lem}

\begin{proof}
Newton's binomial formula yields the identity. For
 $a$ and $a+h$ in $[0,1]$ and $L\ge2$
\[
B(a,h) \le h^L + \sum_{i=2}^{L-1} {L \choose i} a h^2
= h^L + \left(2^L - L -2\right) a h^2\,.
\]
A convexity argument yields $B(a,h)\ge 0$.
\end{proof}

We define the functions $G^N$ mapping $v$ in $c_0$ to 
$G^N(v)$ in $c^0_0$ given by
\begin{equation}
\label{defG}
G^N = F^N - F = F^N_+ - F_+\,,
\quad
G^N(v)(k)
=\alpha A^N(v(k-1)) -  \alpha  A^N(v(k))\,,
\quad
k\ge1\,,
\end{equation}
and $\lin$ and $H$ mapping $(v,x)$ in $c_0 \times c_0^0$
to $\lin(v)x$ and $H(v,x)$ in $c_0^0$ given by
\begin{eqnarray}
\lin(v)x (k) &=& 
\alpha L v(k-1)^{L-1} x(k-1)
-(\alpha L v(k)^{L-1} + \beta) x(k)
+  \beta x(k+1)\,,\quad 
k\ge1\,,\qquad
\nonumber \\
\label{defH}
H(v,x)(k) &=& \alpha  B(v(k-1),x(k-1)) - \alpha B(v(k),x(k))\,,\quad 
k\ge1\,.
\end{eqnarray} 
For $v$ and $v+x$ in $\mathcal{V}$ 
we may use the bounds in Lemmas~\ref{fact} and \ref{lin}. We have
\begin{equation}
\label{diff}
F(v+x) - F(v)  
=F_+(v+x) - F_+(v) + F_-(x)   
= \lin(v) x + H(v,x)\,.
\end{equation}

We derive a limit equation for the fluctuations from 
(\ref{zeq}) and (\ref{bracket}) using (\ref{defG}), (\ref{diff}), and
Lemmas~\ref{fact} and \ref{lin}.
Let $u$ solve (\ref{is}) in $\mathcal{V}$ and $(M(k))_{k\in\NN}$ be independent
real continuous centered Gaussian martingales, determined in law by their
deterministic Doob-Meyer brackets given by
\[
 \langle M(k) \rangle_t = \int_0^t
\left\{ F^{\vphantom{N}}_+(u_s)(k) + F^{\vphantom{N}}_-(u_s)(k) \right\}ds\,. 
\]
The processes $M = (M(k))_{k\ge0}$ and
$\langle M \rangle = \left(\langle M(k) \rangle\right)_{k\in\NN}$
have sample paths with values in $c^0_0$, and
$\lin(u_t) : z \mapsto \lin(u_t)z$
are linear operators on $c^0_0$.
The natural limit equation for the fluctuations is 
the inhomogeneous affine SDE given for $t\ge0$ by
\[
Z_t = Z_0 + \int_0^t \lin(u_s) Z_s \,ds + M_t\,.
\]
We set $\lint = \lin(\ut)$.
For $u_0 = \ut$, (\ref{F}) and $F_+(\ut)=F_-(\ut)$
yield the formulation in Section~\ref{oueq}. 

\section{Main properties of the Ornstein-Uhlenbeck process}
\setcounter{equation}{0}
\label{propOU}

\subsection{Proof of Theorem~\ref{Dbdd}}

Considering (\ref{defDst}) and convexity bounds we have 
\begin{eqnarray*}
\Vert \lint z \Vert_{L_2(w)}^2
&=& \beta^2 \sum_{k\ge1}\left(
L \rho^{L^{k-1}} z(k-1)
-(L \rho^{L^{k}} + 1) z(k)
+  z(k+1)
\right)^2 w(k)^{-1}
\\
&\le&
\beta^2 (2L+2)\biggl(
L\sum_{k\ge1}
\rho^{2L^{k-1}} z(k-1)^2 w(k)^{-1}
+
L\sum_{k\ge1}
\rho^{2L^{k}} z(k)^2 w(k)^{-1}
\\
&&\kern30mm{}
+
\sum_{k\ge1} z(k)^2 w(k)^{-1}
+
\sum_{k\ge1} z(k+1)^2 w(k)^{-1}
\biggr)
\\
&\le&
\beta^2 (2L+2)\biggl(
L d \sum_{k\ge2} z(k-1)^2 w(k-1)^{-1}
+
(L \rho^{2L} +1) \sum_{k\ge1} z(k)^2 w(k)^{-1}
\\
&&\kern30mm{}
+
c^{-1}\sum_{k\ge1} z(k+1)^2 w(k+1)^{-1}
\biggr)
\\
&\le&
\beta^2 (2L+2) \left(L\rho^{2L} + Ld +c^{-1} +1\right) \Vert z \Vert_{L_2(w)}^2\,.
\end{eqnarray*}
The Gronwall Lemma yields uniqueness. For $k\ge1$ we have
\[
(L\rho^L)^{-1}  \pi(k+1) 
\le
\pi(k) =  (L\rho^{L^k})^{-1} \pi(k+1)  \le L^{-1} \rho^L \rho^{-2L^k} \pi(k+1)\,,
\]
\[
\theta^{-1}\theta^{k+1} 
\le  \theta^k 
\le  \theta^{-1} \rho^L  \rho^{-2L^k} \theta^{k+1} \,.
\]
When $B$ is an Hilbertian Brownian motion,
the formula for $Z$ is well-defined and solves the equation.

\subsection{A related birth and death process, and the spectral decomposition}
\label{spectdec}

Considering (\ref{matrixK}), $\mathcal{A} = \lint^*$ is  
the infinitesimal generator   
of the sub-Markovian birth and death process on  the irreducible class $(1,2,\ldots)$
with birth rates $\lambda_k =\beta L \rho^{L^k}$ 
and death rates $\mu_k=\beta$
for $k\ge1$ (killed at rate $\mu_1=\beta$ at state $1$).
The process is well-defined since the rates are bounded.

Karlin and McGregor~\cite{Karlin:57a,Karlin:57b} give a spectral 
decomposition for such  processes, used by
Callaert and Keilson~\cite{Callaert:73a,Callaert:73b} 
and van Doorn~\cite{Doorn:85}
to study exponential ergodicity
properties.
The state space in these works is $(0,1,2,\dots)$,
possibly extended by an absorbing barrier 
or graveyard state at
$-1$. We consider  $(1,2,\ldots)$ and adapt their notations to this simple shift.

The potential coefficients (\cite{Karlin:57a} eq.~(2.2), \cite{Doorn:85} eq.~(2.10)) are
given  by 
\[
\pi(k) 
= {\lambda_1 \cdots \lambda_{k-1} \over \mu_2 \cdots \mu_k}
=  L \rho^{L^1}\cdots  L \rho^{L^{k-1}} =   L^{k-1} \rho^{(L^k-L)/(L-1)},
\qquad k\ge1\,,
\]
and solve the detailed balance
equations $\mu_{k+1} \pi(k+1) = \lambda_{k} \pi(k)$ with
$\pi(1)=1$. 

The equation ${\mathcal A}Q(x) = - x Q(x)$
for an eigenvector $Q(x) = (Q_n(x))_{n\ge1}$ of eigenvalue $-x$ 
yields $\lambda_1 Q_2(x) = (\lambda_1 + \mu_1 -x) Q_1(x)$ and
$\lambda_n Q_{n+1}(x) = (\lambda_n + \mu_n -x) Q_n(x) - \mu_n Q_{n-1}(x)$ for $n\ge2$.
With the natural convention $Q_0=0$
and choice $Q_1=1$, we obtain
inductively  $Q_n$ as the  polynomial of degree $n-1$ satisfying
\[
- x  Q_n(x)
=  \beta  Q_{n-1}(x)  -  \left(\beta L \rho^{L^n}   + \beta  \right) Q_n(x) 
+ \beta  L \rho^{L^n}  Q_{n+1}(x)\,,
\qquad  n\ge1\,.
\]
These recursions correspond to 
\cite{Karlin:57a}~eq.~(2.1)
and \cite{Doorn:85}~eq.~(2.15). As stated there, such 
a sequence of polynomials is orthogonal with respect to a probability measure
$\psi$ on $\RR_+$ and
\[
\int_0^\infty Q_i(x)^2\, \psi(dx) = \pi(i)^{-1}\,, 
\quad
\int_0^\infty Q_i(x) Q_j(x)\, \psi(dx) = 0\,, \qquad 
i,j\ge1\,,\ i\neq j\,,
\]
or in matrix notation
$\int_0^\infty Q(x) Q(x)^*\, \psi(dx) = \mathrm{diag}(\pi^{-1})$.

Let $P_t = (p_t(i,j))_{i,j\ge1}$ denote the sub-stochastic 
transition matrix for $\mathcal{A}$. The adjoint matrix
$P_t^*$ 
is the fundamental solution for the forward Kolmogorov equation 
$\dot z_t = \mathcal{A}^* z_t = \lint z_t$.
The representation formula of 
Karlin and McGregor~\cite{Karlin:57a,Karlin:57b} 
(see (1.2) and (2.18) in \cite{Doorn:85})
yields
\begin{equation}
\label{repr}
\mathrm{e}^{\lint t} = P_t^* = (p^*_t(i,j))_{i,j\ge1}\,,
\qquad
p^*_t(i,j) = p_t(j,i) = 
\pi(i) \int_0^\infty 
\mathrm{e}^{-xt} Q_i(x) Q_j(x)\, \psi(dx)\,,
\end{equation}
or in  matrix notation
$\mathrm{e}^{\lint t} 
= \mathrm{diag}(\pi) \int_0^\infty  \mathrm{e}^{-xt} Q(x) Q(x)^* \,\psi(dx)$.

The probability measure $\psi$ is called the spectral measure, its support $S$
is called the spectrum, and we set  $\gamma=\min S$. The 
Ornstein-Uhlenbeck process in Theorem~\ref{Dbdd}~(b) and its invariant law
and its covariance matrix in Theorems~\ref{ou.erg} and \ref{eqfclt} 
can be written
\begin{eqnarray}
\label{explou}
Z_t &=&
\mathrm{diag}(\pi)
\int_S  \mathrm{e}^{-xt}   Q(x)^*  \left(Z_0
+\int_0^t\mathrm{e}^{xs} \,dB_s 
\right) Q(x) \,\psi(dx)\,,
\\
\label{limlaw}
\int_0^\infty \mathrm{e}^{\lint t}\,dB_t
&=&
\mathrm{diag}(\pi) \int_S 
\left(Q(x)^*\int_{0}^\infty  \mathrm{e}^{-xt} \,dB_t \right) Q(x) \,\psi(dx)\,,
\\
\label{cov}
\int_0^\infty \mathrm{e}^{\lint t} \mathrm{diag}(\tilde v) \mathrm{e}^{\lint^* t}\,dt
&=&
\mathrm{diag}(\pi)
\int_{S^2} 
{Q(x)^* \mathrm{diag}(\tilde v) Q(y) \over x+y}\,
Q(x) Q(y)^*  \,\psi(dx)\psi(dy)\,
\mathrm{diag}(\pi).\qquad
\end{eqnarray}

\subsection{The spectral gap, exponential stability, and ergodicity}
\label{prf.spe.gap}

{\em Proof of Theorem~\ref{spe.gap}\/}.
The potential coefficients $(\pi(k))_{k\ge1}$ solve the detailed balance equations 
for $\mathcal{A}$ and hence 
$\lint = \mathcal{A}^*$ is self-adjoint in 
$L_2(\pi)$.

For the spectral gap,
we follow Van Doorn~\cite{Doorn:85}, Section~2.3.
The orthogonality properties imply that for $n\ge1$,
 $Q_n$ has $n-1$ distinct zeros $0 < x_{n,1}<\ldots<x_{n,n-1}$ such that
$x_{n+1,i} <x_{n,i} < x_{n+1,i+1}$
for $1\le i \le n-1$. Hence $\xi_i = \lim_{n\to\infty }x_{n,i} \ge 0$
exists, $\xi_i\le \xi_{i+1}$, and $\sigma = \lim_{i\to\infty}\xi_i$ exists
in $[0,\infty]$. 
Theorem 5.1 in \cite{Doorn:85} establishes that
$\gamma>0$ if and only if $\sigma>0$,
Theorem~5.3~(i) in \cite{Doorn:85} that
$\sigma = \beta>0$,  and
Theorem 3.3 in \cite{Doorn:85}  that $\gamma = \xi_1 \le \sigma$.
(Estimating $\xi_1$ is impractical.)

For the exponential stability, we have
$\Vert z_t \Vert_{L_2(\pi)}^2 = \left(\mathrm{e}^{\lint t} z_0, 
\mathrm{e}^{\lint t} z_0 \right)_{L_2(\pi)}$.
The fact that
$\mathrm{e}^{\lint t}$ is self-adjoint in $L_2(\pi)$ and 
the spectral representation (\ref{repr}) yield
\begin{eqnarray*}
\left(\mathrm{e}^{\lint t} z_0, 
\mathrm{e}^{\lint t} z_0 
\right)_{L_2(\pi)}
&=&
\left( z_0, 
\mathrm{e}^{ 2\lint t} z_0 
\right)_{L_2(\pi)}
=
\int_S \mathrm{e}^{- 2xt} z_0^* Q(x) Q(x)^* z_0\,\psi(dx) 
\\
&\le& \mathrm{e}^{- 2\gamma t}  
\int_S z_0^* Q(x) Q(x)^* z_0\,\psi(dx) 
=
\mathrm{e}^{- 2\gamma t} \left( z_0,z_0 \right)_{L_2(\pi)}.
\end{eqnarray*}
We refer to 
Callaert and Keilson~\cite{Callaert:73b} Section 10 for related results.

\medskip
\noindent
{\em Proof of Theorem~\ref{linis.exp} (non self-adjoint case)\/}.
It is similar to and simpler than the proof for Theorem~\ref{glexst}
in the interactive case $L\ge2$, and we wait till that point to give it.

\medskip
\noindent{\em Proof of Theorem~\ref{ou.erg}\/}.
We use the uniqueness result and
explicit formula for $Z$ in Theorem~\ref{Dbdd},
and Theorem~\ref{spe.gap} or \ref{linis.exp}.

\section{Exponential stability for the nonlinear system}
\setcounter{equation}{0}
\label{expstab}

\subsection{Some comparison results}

Considering  (\ref{diff}), $\lint = \lin(\ut)$  and $F(\ut)=0$,
if $u$ is a solution of (\ref{is}) in $\mathcal{V}$
starting at $u_0$
then $y = u-\ut$ is a solution to the 
recentered equation starting at $y_0=u_0-\ut$ given by 
\begin{eqnarray}
\label{center}
\kern-5mm
\dot y_t (k) 
&=&
\lint y_t(k) + H(\ut, y_t)(k)
\nonumber\\
&=&
\beta L \rho^{L^{k-1}}y_t(k-1) + \alpha  B(\ut(k-1),y_t(k-1))  
\nonumber\\
&&\kern3mm{}
-\left(\beta L \rho^{L^k}y_t(k) + \alpha B(\ut(k),y_t(k)) + \beta y_t(k) \right) 
+  \beta y_t(k+1)\,,
\qquad k\ge1\,,
\end{eqnarray}
and if $u_0$ is in $\mathcal{V}\cap \ell_1$ then $u$ is in $\mathcal{V}\cap \ell_1$
and hence $y$ is in $\ell_1^0$ and for $k\ge1$ 
\begin{equation}
\label{sum.cen}
\dot y_t(k) + \dot y_t(k+1) + \cdots\,
=  \beta L \rho^{L^{k-1}} y_t(k-1) + \alpha B(\ut(k-1), y_t(k-1))  - \beta y_t(k)\,.
\end{equation}
Reciprocally, if $y$ 
is a solution to the  recentered  equation (\ref{center}) starting at 
$y_0$ such that $y_0+ \ut$ is in $\mathcal{V}$, then $u = y+ \ut$ 
is a solution of (\ref{is}) in $\mathcal{V}$
starting at $u_0=y_0+ \ut$. Then $-\ut \le y \le 1-\ut$ and $-1 < y  < 1$.
For $y_0+ \ut$ in $\mathcal{V}\cap\ell_1$ we have  $y$ in $\ell_1^0$. 

\begin{lem}
\label{order}
Let $u$ and $v$ be two solutions for (\ref{is}) 
in $\mathcal{V}$ such that $u_0 \le v_0$. Then $u_t \le v_t$ for $t\ge0$.
Let $y_0+ \ut$ be in $\mathcal{V}$ and $y$ solve (\ref{center}).
If $y_0\ge0$ then $y_t\ge0$ and if $y_0\le0$ then $y_t\le0$ for $t\ge0$.
\end{lem}

\begin{proof}
Lemma 6  in \cite{Vved:96} yields the result for  (\ref{is}) (the proof written 
for $L=2$ is valid for $L\ge1$).
The result for (\ref{center}) follows by consideration of
the  solutions $u=y+\ut$ and $\ut$ for (\ref{is}).
\end{proof}

We shall compare solutions of the nonlinear equation (\ref{center}) 
and of certain linear equations.

\begin{lem}
\label{kolm}
Let $\hat\mathcal{A}$ be the generator of the 
sub-Markovian birth and death process
with birth rate $\hat\lambda_k \ge 0$ 
and death rate $\beta$ at $k\ge1$. Let $\sup_k \hat\lambda_k<\infty$.
In $\ell_1^0$ the linear operator 
\[
\hat\mathcal{A}^*x(k) = 
\hat\lambda_{k-1} x(k-1)
-(\hat\lambda_{k} + \beta ) x(k)
+  \beta x(k+1)\,,
\qquad
k\ge1\,,
\]
is bounded and there exists a unique $z=(z_t)_{t\ge0}$ given by 
$z_t = \mathrm{e}^{\hat\mathcal{A}^* t} z_0$
solving
the forward Kolmogorov equation $\dot z=\hat\mathcal{A}^*z$.
If $z_0\ge0$ then $z_t\ge0$ and if
$z_0\le0$ then $z_t\le0$.
For $k\ge1$,
$
\dot z_t(k) + \dot z_t(k+1) + \cdots\,
=   \hat\lambda_{k-1} z_t(k-1)  - \beta z_t(k)
$.
\end{lem}

\begin{proof}
The  operator norm in $\ell_1^0$ of $\hat\mathcal{A}^*$ 
is bounded by $2(\sup_k \hat\lambda_k + \beta)$, hence existence and uniqueness.
Uniqueness and linearity imply that if $z_0=0$ then $z_t=0$
and else if
$z_0\ge0$ then 
$z_t \Vert z_0\Vert_1^{-1}$ is the instantaneous
law of the process starting at  $z_0\Vert z_0\Vert_1^{-1}$
and hence $z_t\ge0$. If $z_0\le0$ then $- z$ solves the equation
starting at $- z_0 \ge 0$ and hence $- z_t \ge 0$.
\end{proof}

\begin{lem}
\label{cen.lin}
Let $L\ge2$ and $y=(y_t)_{t\ge0}$ solve (\ref{center}) with  $y_0+ \ut$ in 
$\mathcal{V}\cap\ell_1$.
Under the assumptions  of Lemma~\ref{kolm},
let $z=(z_t)_{t\ge0}$  solve 
$\dot z  = \hat\mathcal{A}^*z $ with $z_0$ in $\ell_1^0$ and
$h=(h_t)_{t\ge0}$ be given by
\[
h = (h(k))_{k\ge1}\,,
\qquad
h(k) = z(k) + z(k+1) + \cdots\,
- (y(k) + y(k+1) + \cdots\,).
\]

\noindent
(a) Let $\hat\lambda_k \ge \beta L \rho^{L^k} 
+  \alpha\!\left(1 + \left(2^L-L-2\right) \ut(k)\right)$  for $k\ge1$,
$y_0\ge0$, and $h_0\ge0$. Then $h_t\ge0$ for $t\ge0$.

\noindent
(b) Let $\hat\lambda_k \ge \beta L \rho^{L^k}$  for $k\ge1$,
$y_0\le0$, and $h_0\le0$.   Then $h_t\le0$ for $t\ge0$.
\end{lem}

\begin{proof}
We prove (a).
For $\varepsilon >0$ let $\hat\mathcal{A}^*_\varepsilon$ correspond to
$\hat\lambda_k^\varepsilon = \hat\lambda_k + \varepsilon$.
The operator norm in $\ell_1^0$ of
$\hat\mathcal{A}^*_\varepsilon -\hat\mathcal{A}^*$ is bounded by
$2 \varepsilon$, hence
$\lim_{\varepsilon\to0} \mathrm{e}^{\hat\mathcal{A}_\varepsilon^* t} z_0 =z_t$ 
in $\ell_1^0$ and we may 
assume that 
$\hat\lambda_k > \beta L \rho^{L^k}
+  \alpha\!\left(1 + \left(2^L-L-2\right) \ut(k)\right)$ for $k\ge1$.
Since
$z_t = \mathrm{e}^{\hat\mathcal{A}^* t} z_0$ depends continuously on $z_0$ 
in $\ell_1^0$
we may assume $h_0 > 0$. 

Let $\tau = \inf\{t\ge0 : \{k\ge1 : h_{t}(k)=0\} \neq\emptyset \}$  be
the first time when $h(k)=0$ for some $k\ge1$. 
Then $\tau>0$ and 
if $\tau =\infty$ the proof is ended. Else, 
Lemma \ref{kolm} and (\ref{sum.cen}) yield
\begin{eqnarray*}
\dot h_\tau(k) &=&
\hat\lambda_{k-1}y_\tau(k-1)  - \beta L \rho^{L^{k-1}} y_\tau(k-1)
- \alpha B(\ut(k-1), y_\tau(k-1))
\nonumber\\
&&\kern5mm{}
+\hat\lambda_{k-1} (z_\tau(k-1) - y_\tau(k-1) )
-\beta  (z_\tau(k) -  y_\tau(k))\,.
\end{eqnarray*}
Lemma~\ref{order} yields $y \ge 0$ and 
Lemma~\ref{lin}  and $y\le1$ yield
\begin{eqnarray*}
B(\ut(k-1), y(k-1))&\le&
y(k-1)^L + \left(2^L -L -2\right)\ut(k-1) y(k-1)^2
\\
&\le& \left(1 + \left(2^L -L -2\right)\ut(k-1)\right) y(k-1)\,,
\end{eqnarray*}
hence
$
\hat\lambda_{k-1}y(k-1)  - \beta L \rho^{L^{k-1}} y(k-1)
- \alpha B(\ut(k-1), y(k-1)) \ge0
$
with equality only when $y(k-1) = 0$.
For $k$ in $\mathcal{K} = \{k\ge1 : h_\tau(k)=0 \} \neq \emptyset$ we have
\[
z_\tau(k-1) - y_\tau(k-1) =  h_\tau(k-1)\ge 0\,,
\quad
z_\tau(k) - y_\tau(k)  = - h_\tau(k+1) \le 0 \,,
\]
with equality if only if 
$k-1$ is in $\mathcal{K}\cup\{0\}$ and $k+1$ is in $\mathcal{K}$.
Hence $\dot h_\tau(k) \ge 0$.
Moreover $h_t(k) > 0$ for $t<\tau$
and $h_\tau(k) = 0$ imply $\dot h_\tau(k) \le 0$, hence
$\dot h_\tau(k) = 0$, and the above signs and equality cases  yield that
$z_\tau(k-1)=y_\tau(k-1)=0$ and
$k-1$ is in $\mathcal{K}\cup\{0\}$ and $k+1$ is in $\mathcal{K}$.
By induction $z_\tau(i) = y_\tau(i)=0$ for $i\ge1$ which implies 
$z_{t} = y_{t}=0$ for $t\ge \tau$.

The proof for (b) is similar and involves obvious changes of sign.
We may assume $\hat\lambda_k > \beta L \rho^{L^k}$ which suffices to conclude since
Lemma~\ref{lin}  yields $B(\ut(k-1), y(k-1))\ge0$.
\end{proof}

\begin{lem}
\label{equiv}
For any $0<\theta<1$ there exists $K_\theta<\infty$ such that
for $x$ in $L_2(g_\theta)\subset \ell^0_1$ 
\[
\left\Vert (x(k)+x(k+1)+\cdots)_{k\ge1} \right \Vert_{L_2(g_\theta)}
\le
 K_\theta \Vert x \Vert_{L_2(g_\theta)}\,.
\]
\end{lem}

\begin{proof}
Using a classical convexity inequality
\begin{eqnarray*}
&&\kern-2mm
\sum_{k\ge1}  (x(k)+ x(k+1) + \cdots\, )^2 \theta^{-k} 
\\
&&\kern2mm {} \le
\sum_{k\ge1} n\!\left(
x(k)^2 + x(k+1)^2 + \cdots + x(k+n-2)^2 
+ (x(k+n-1)+ x(k+n) + \cdots\, )^2
\right) \theta^{-k} 
\\
&&\kern2mm {} \le 
n \!\left(1+\theta + \cdots+ \theta^{n-2} \right)\sum_{k\ge1}  x(k)^2 \theta^{-k}
+ n\,  \theta^{n-1}\sum_{k\ge1}  (x(k)+ x(k+1) + \cdots\, )^2 \theta^{-k}.
\end{eqnarray*} 
We take $n$ large enough that $n\theta^{n-1}<1$ and
$
K_\theta^2 = (1 - n \theta^{n-1})^{-1} n (1 - \theta^{n-1})(1 - \theta)^{-1}\,.
$
\end{proof}

\subsection{Proofs of Theorems~\ref{glexst} and \ref{linis.exp}}

{\em Proof of Theorem~\ref{glexst} for $L\ge2$\/}.
Let $u_0$ be in $\mathcal{V} \cap L_2(g_\theta)$. Then 
$u^-_0 = \min\{u_0,\ut\}$ and $u^+_0 = \max\{u_0,\ut\}$
are in $\mathcal{V} \cap L_2(g_\theta)$.
Theorem~\ref{isw} yields that the 
corresponding solutions $u^-$ and $u^+$ for (\ref{is})
are in $\mathcal{V}\cap L_2(g_\theta)$.
Lemma~\ref{order} yields that $u^-_t  \le u_t \le u^+_t$ and 
$u^-_t  \le \ut \le u^+_t$ 
for $t\ge0$. Then
\[
y= u-\ut\,,
\qquad
y^+ = u^+ - \ut \ge 0\,,
\qquad
y^- =  u^- - \ut \le 0\,,
\]
solve (\ref{center}), and termwise
\begin{equation}
\label{trmin}
|y_0| = \max\{y^+_0 , - y^-_0\}\,,
\qquad
|y_t| \le \max\{y^+_t , - y^-_t\} \,,
\quad t\ge0\,.
\end{equation}

We consider the birth and death process with generator $\hat\mathcal{A}$
defined in Lemma~\ref{kolm}
with 
\[
\hat \lambda_k 
= \max\left\{ 
\beta L \rho^{L^{k}} +   \alpha\!\left(1 + \left(2^L-L-2\right) \ut(k)\right) , 
\beta \theta 
\right\}\,,
\qquad
k\ge1\,,
\]
which satisfies the assumptions of Lemma~\ref{cen.lin} (a) and (b).
We perform 
the same spectral study 
as in Sections~\ref{spectdec} and \ref{prf.spe.gap}, all notions being similar
and denoted using a hat.

For $\rho \le \theta <1$ we have $\alpha \le \beta \theta$ and hence
$\hat \lambda_k$ is equivalent to $\beta \theta$ for large $k$, hence
Theorem~5.3~(i) in \cite{Doorn:85} yields that
$0<\hat\gamma  \le \hat\sigma 
= \left(\sqrt{\beta} - \sqrt{\beta\theta} \right)^2 
= \beta \left(1 - \sqrt{\vphantom{\beta}\theta} \right)^2$, and
moreover 
\[
\theta^{k-1} 
\le
\hat\pi(k) = \theta^{k-1} 
\prod_{i=1}^{k-1}
\max\left\{ 
\theta^{-1}L \rho^{L^{k}} 
+ \theta^{-1}\rho\!\left(1 + \left(2^L-L-2\right) \ut(k)\right),
1 \right\} 
\]
and the product converges using simple criteria. Hence
$\hat\pi(k) = O(\theta^{k})$ and $\theta^{k} = O(\hat\pi(k))$ and
Lemma~\ref{eqiv.l2w} yields that there exists
$c>0$ and $d>0$ such that 
$c^{-1} \Vert \cdot \Vert_{L_2(\hat\pi)} \le 
\Vert \cdot \Vert_{L_2(g_\theta)} \le d \Vert \cdot \Vert_{L_2(\hat\pi)}$.
The version of Theorem~\ref{spe.gap} for the the above process yields that
if $z$ solves $z=\hat\mathcal{A}^*z$ in  $ L_2(g_\theta)$ then
\[
\Vert z_t \Vert_{L_2(g_\theta)} 
\le d \Vert z_t \Vert_{L_2(\hat\pi)}
\le \mathrm{e}^{-\hat\gamma t}  d  \Vert z_0 \Vert_{L_2(\hat\pi)}
\le  \mathrm{e}^{-\hat\gamma t} cd \Vert z_0 \Vert_{L_2(g_\theta)}\,.
\]
Hence if $z^+$ solves $z^+= \hat\mathcal{A}^* z^+$ starting at $z_0^+ =y_0^+ \ge 0$
then Lemma~\ref{cen.lin} (a) and Lemma~\ref{equiv} yield
\begin{eqnarray*}
\Vert y^+_t \Vert_{L_2(g_\theta)} 
&\le &
\Vert (y^+_t(k) + y^+_t(k+1) +\cdots\, )_{k\ge1} \Vert_{L_2(g_\theta)}
\\
&\le &
\Vert (z^+_t(k) +  z^+_t(k+1) +\cdots\, )_{k\ge1} \Vert_{L_2(g_\theta)}
\\
&\le &
K_\theta \Vert z^+_t \Vert_{L_2(g_\theta)} 
\le \mathrm{e}^{-\hat\gamma t} cd K_\theta \Vert y^+_0 \Vert_{L_2(g_\theta)}\,,
\end{eqnarray*}
and similarly 
if $z^-$ solves $z^-= \hat\mathcal{A}^* z^-$ starting at $z_0^- =y_0^- \le0$
then Lemma~\ref{cen.lin} (b) and Lemma~\ref{equiv} yield
$
\Vert y^-_t \Vert_{L_2(g_\theta)} \le  
\mathrm{e}^{-\hat\gamma t} cd K_\theta \Vert y^-_0 \Vert_{L_2(g_\theta)}
$.
We set $\gamma_\theta = \hat\gamma$ and $C_\theta = cd K_\theta$.
Considering (\ref{trmin}),
\[
\Vert y_t \Vert_{L_2(g_\theta)}^2 
\le
\Vert y^+_t \Vert_{L_2(g_\theta)}^2 + \Vert y^-_t \Vert_{L_2(g_\theta)}^2 
\le \mathrm{e}^{- 2\gamma_\theta t} C_\theta^2
\left(
\Vert y^+_0 \Vert_{L_2(g_\theta)}^2 + \Vert y^-_0 \Vert_{L_2(g_\theta)}^2 
\right)
\]
and 
we complete the proof by remarking that
for $k\ge1$, either $y^+_0(k) = y_0(k)$ and  $y^-_0(k) = 0$ or
$y^-_0(k) = y_0(k)$ and $y^+_0(k) = 0$, and hence
$
\Vert y^+_0 \Vert_{L_2(g_\theta)}^2 + \Vert y^-_0 \Vert_{L_2(g_\theta)}^2 
 = 
\Vert y_0 \Vert_{L_2(g_\theta)}^2$.
\medskip

\noindent
{\em Proof of Theorem~\ref{linis.exp} and of Theorem~\ref{glexst} for  $L=1$\/}.
The linearization (\ref{linis})
of Equation (\ref{is}) is obtained from
Equation (\ref{center}) by replacing the nonlinear functions $B$ and $H$ by $0$,
and coincides with (\ref{center}) for $L=1$. Likewise,
the equation for (\ref{linis}) corresponding to (\ref{sum.cen})
is obtained by omitting the term $\alpha B(\ut(k-1), y_t(k-1))$.
We obtain a result for the linear equation (\ref{linis}) 
corresponding to Lemma~\ref{cen.lin} (a) and (b)
under the sole assumption
$\hat \lambda_k \ge \beta L \rho^{L^{k}}$ for $k\ge1$.
The proof proceeds as for Theorem~\ref{glexst} for $L\ge2$ 
with the difference that
$\hat \lambda_k = \max\left\{\beta L \rho^{L^{\smash{k}}}, \beta \theta \right \}$.
We have $\hat \lambda_k$  equal to $\beta \theta$ for large $k$
for $0<\theta< 1$ when  $L\ge2$ and for $\rho \le \theta <1$
when $L=1$.

\section{Tightness estimates and the functional central limit theorem}
\label{cfp}
\setcounter{equation}{0}

\subsection{Finite horizon bounds for the process: proof of Lemma~\ref{l2bds}}
\label{sprl2bds}

We use Lemma \ref{dyn}. Considering (\ref{zeq}) and (\ref{defG}), 
\begin{equation}
\label{znavecg}
Z^N_t 
= Z^N_0 + M^N_t + N^{1/2} \int_0^t G^N(R^N_s)\,ds
+\int_0^t 
N^{1/2}\left(F (R^N_s) - F(\ut) \right)  ds
\end{equation}
where Lemma~\ref{fact} yields that 
\[
G^N(R^N_s)(k) =
\alpha \!\left(
A^N \!\left(R^N_s (k-1)\right) 
- A^N \!\left(R^N_s(k) \right) 
\right)
= N^{-1} O\left( R^N_s(k-1)  +  R^N_s(k) \right)
\]
and hence for some $K<\infty$
\begin{equation}
\label{grz}
\left\Vert G^N(R^N_s)\right\Vert_{L_2(g_\theta)}
\le
N^{-1} K \left\Vert R^N_s \right\Vert_{L_2(g_\theta)}
\end{equation}
where
\begin{equation}
\label{decr1}
\left\Vert R^N_s \right\Vert_{L_2(g_\theta)}
\le
\left\Vert \ut \right\Vert_{L_2(g_\theta)} + 
N^{-1/2}  \left\Vert Z^N_s\right\Vert_{L_2(g_\theta)}.
\end{equation}
The mapping $F$ being Lipschitz (Theorem~\ref{isw}), the Gronwall Lemma 
yields that for some $K_T<\infty$
\[
\sup_{0\le t\le T}\left\Vert Z^N_t \right\Vert_{L_2(g_\theta)}
\le
K_T\biggl(\left\Vert Z^N_0 \right\Vert_{L_2(g_\theta)} 
+ \sup_{0\le t\le T}\left\Vert M^N_t \right\Vert_{L_2(g_\theta)}
+ N^{-1/2} \left\Vert \ut \right\Vert_{L_2(g_\theta)}
\biggr).
\]
We conclude using the Doob inequality, (\ref{bracket}), (\ref{defG}),
the bounds (\ref{grz}) and (\ref{decr1}), and (see Theorem~\ref{isw})
\begin{equation}
\label{brkbd}
\left\Vert F_+(R^N_s) + F_-(R^N_s) \right\Vert_{L_2(g_\theta)}
\le
K\left\Vert R^N_s \right\Vert_{L_2(g_\theta)}.
\end{equation}

\subsection{Infinite horizon bounds for the marginals: proof of Lemma~\ref{l2bdsinf}}

Let  $U_h (v)$ be the solution of (\ref{is}) at time 
$h\ge0$ with initial value $v$ in $\mathcal{V}$, in particular $\ut= U_h(\ut)$, and
$Z^N_{t_0,h} = N^{1/2}\left(R^N_{t_0+h}-  U_h(R^N_{t_0}) \right)$ for $t_0\ge0$.
We have
$
Z^N_{t_0+h}=  Z^N_{t_0,h} + N^{1/2}\left( U_h(R^N_{t_0}) - \ut \right)
$
and Theorem~\ref{glexst} yields that
\begin{equation}
\label{decineg}
\left\Vert Z^N_{t_0+h}  \right\Vert_{L_2(g_\theta)}
\le 
\left\Vert
Z^N_{t_0,h}
\right\Vert_{L_2(g_\theta)}
+
\mathrm{e}^{- \gamma_\theta h} C_\theta 
\left\Vert
Z^N_{t_0}
\right\Vert_{L_2(g_\theta)}.
\end{equation}

The conditional law of
$(Z^N_{t_0,h})_{h\ge0}$ 
given $R^N_{t_0} =r$
is the law of $Z^N$
started with $R^N_0 = u_0 =r$,
the empirical fluctuation process centered on $U(r)$ and starting at $0$.
We reason as in Section~\ref{sprl2bds}, using additionally
(\ref{decineg}) on the bound  (\ref{decr1}) with $s=t_0+h$.
We obtain that for some $K_T<\infty$
\[
\sup_{0\le h \le T} \left\Vert Z^N_{t_0,h} \right\Vert_{L_2(g_\theta)} 
\le K_T  \biggl(
N^{-1}  C_\theta 
\left\Vert
Z^N_{t_0}
\right\Vert_{L_2(g_\theta)}
+
\sup_{0\le h \le T} \left\Vert M^N_{t_0+h} - M^N_{t_0} \right\Vert_{L_2(g_\theta)}
+ N^{-1/2} \left\Vert \ut \right\Vert_{L_2(g_\theta)}
\biggr)
\]
and then that for some $L_T<\infty$ we have for $0\le h \le T$
\begin{equation}
\label{bonineg}
\EE\left(
\left\Vert Z^N_{t_0+h}  \right\Vert_{L_2(g_\theta)}^2
\right)
\le 
L_T
+
2 ( K_T N^{-1} + \mathrm{e}^{- \gamma_\theta h})^2 C_\theta^2 
\,\EE\left(
\left\Vert
Z^N_{t_0}
\right\Vert_{L_2(g_\theta)}^2
\right).
\end{equation}

We fix $T$ large enough for 
$8\mathrm{e}^{- 2\gamma_\theta T}C_\theta^2  \le \varepsilon < 1$.
Uniformly for $N \ge K_T \mathrm{e}^{\gamma_\theta T}$, 
for $m\in\NN$
\[
\EE\left(
\left\Vert Z^N_{\smash{(m+1)}T}  \right\Vert_{L_2(g_\theta)}^2
\right)
\le 
L_T
+
\varepsilon
\,\EE\left(
\left\Vert
Z^N_{mT}
\right\Vert_{L_2(g_\theta)}^2
\right)
\]
and  by induction
\[
\EE\left(
\left\Vert Z^N_{mT}  \right\Vert_{L_2(g_\theta)}^2
\right)
\le 
L_T
\sum_{j=1}^m
\varepsilon^{j-1}
+
\varepsilon^m
\,\EE\left(
\left\Vert
Z^N_{0}
\right\Vert_{L_2(g_\theta)}^2
\right)
\le 
{L_T \over 1 - \varepsilon} + 
\EE\left(\left\Vert
Z^N_{0}
\right\Vert_{L_2(g_\theta)}^2
\right),
\]
and (\ref{bonineg}) yields
\[
\sup_{0 \le h \le T}
\EE\left(
\left\Vert Z^N_{mT+h}  \right\Vert_{L_2(g_\theta)}^2
\right)
\le 
L_T
+
8 C_\theta^2 
\,\EE\left(
\left\Vert
Z^N_{mT}
\right\Vert_{L_2(g_\theta)}^2
\right),
\]
hence 
\[
\sup_{t\ge0}\EE\left(
\left\Vert Z^N_{t}  \right\Vert_{L_2(g_\theta)}^2
\right) 
\le
L_T
+
8 C_\theta^2 
\left(
{L_T \over 1 - \varepsilon} + 
\EE\left(\left\Vert
Z^N_{0}
\right\Vert_{L_2(g_\theta)}^2
\right)
\right).
\]

\subsection{Bounds on the invariant laws: proof of Lemma~\ref{tight.in}}

Ergodicity
and the Fatou Lemma yield that for $Z_\infty^N$ 
distributed according to the invariant law
\[
\EE\left(\left\Vert Z^N_\infty \right\Vert_{L_2(g_\theta)}^2\right)
\le
\liminf_{t\ge0} 
\EE\left(\left\Vert Z^N_t \right\Vert_{L_2(g_\theta)}^2\right)
\le
\sup_{t\ge0} 
\EE\left(\left\Vert Z^N_t \right\Vert_{L_2(g_\theta)}^2\right)
\]
and considering Lemma~\ref{l2bdsinf} the proof will be complete
as soon as we show that we can choose $R^N_0$ in $\mathcal{V}^N$ such
that
\begin{equation}
\label{initgood}
\limsup_{N\ge L} 
\EE\left(\left\Vert Z^N_0 \right\Vert_{L_2(g_\theta)}^2\right) < \infty\,.
\end{equation}

We consider $L\ge2$, the case $L=1$ being similar.
Let $R_0^N= (R_0^N(k))_{k \in \NN}$ with
\[
R_0^N(k) = iN^{-1}
\;\; \mbox{for}\;\;  - (2N)^{-1} < \ut (k) - iN^{-1} \le  (2N)^{-1}\,,
\qquad
i\in\{0,1,\ldots,N\}\,,
\]
and
\[
k(N) = \inf\{k\ge 1 : R_0^N(k)=0  \}= \inf\{k\ge 1 : \ut (k) \le (2N)^{-1}\}\,.
\]
Since for $x\ge 0$ and $0< y\le 1$
\begin{eqnarray*}
y = \rho^{(L^x-1)/(L-1)}
&\Leftrightarrow&
x = \log\left(1 + (L-1)\log y / \log \rho \right) / \log L
\\
&\Leftrightarrow&
\theta^{-x} = \left(1 + (L-1) \log y /\log \rho \right)^{-\log\theta / \log L}
\end{eqnarray*}
we have
$
k(N) = \inf\left\{k\in\NN 
: k \ge \log\left(1 + (L-1)\log \left( (2N)^{-1}\right) / \log \rho \right) / \log L
\right\}$.
Then
\[
 \left\Vert Z^N_0 \right\Vert_{L_2(g_\theta)}^2
=
N \sum_{k=1}^{k(N)-1} \left(R^N_0(k) - \ut (k)\right)^2 \theta^{-k}
+
N \sum_{k\ge k(N)} \ut (k)^2 \theta^{-k},
\]
\[
N \sum_{k=1}^{k(N)-1} \left(R^N_0(k) - \ut (k)\right)^2 \theta^{-k}
\le (4N)^{-1}\,{\theta^{-k(N)}-\theta^{-1}  \over \theta^{-1} -1 }
= O \left( N^{-1} (\log N)^{\smash{-\log \theta / \log L}}  \right),
\]
and for large enough $N$ (and hence $k(N)$)
\begin{eqnarray*}
N \sum_{k\ge k(N)} \ut (k)^2 \theta^{-k}
&=&
N\sum_{k\ge k(N)} \rho^{2(L^k-1)/ (L-1)}\theta^{-k}
\\
&=&
  N 
\rho^{2(L^{k(N)}-1)/ (L-1)}
 \sum_{k\ge k(N)} \rho^{2(L^k-L^{k(N)})/ (L-1)}\theta^{-k}
\\
&\le&
(4N)^{-1}\,\sum_{j\ge0 } \rho^{2L^{k(N)}(L^j-1)/ (L-1)}\theta^{-(j+k(N))}
\\
&\le&
(4N)^{-1}\,\sum_{j\ge0 } \rho^{L^{k(N)}(L^j-1)/ (L-1)} = o(N^{-1}).
\end{eqnarray*}
Hence (\ref{initgood}) holds and the proof is complete.

\subsection{The functional CLT: Proof of Theorem~\ref{eqfclt}}
\label{pclt}

Lemma~\ref{tight.in} and the Markov inequality imply that
in equilibrium
$(Z_0^N)_{N\ge L}$ is asymptotically tight for the weak
topology of $L_2(g_\rho)$, for which all bounded sets are relatively compact.
We consider a subsequence of $N\ge L$.
Let $(N_j)_{j\ge1}$ denote a further 
subsequence such that $(Z_0^{N_j})_{j\ge1}$ converges 
in law to some square-integrable
$Z_0^\infty$
in $L_2(g_\rho)$.
We decompose the rest of the proof in three steps.

\noindent{\em Step 1\/}.
We  prove that $(Z^{N_j})_{j\ge1}$ is tight in 
$\DD(\RR_+, L_2(g_{\rho}))$ with the Skorokhod
topology, where 
$L_2(g_{\rho})$ is considered with its 
non-metrizable
weak topology.
The compact subsets of $L_2(g_{\rho})$ 
are metrizable and hence Polish, a fact yielding
tightness criteria. 
We easily deduce from Theorem~4.6 and 3.1 in 
Jakubowski~\cite{Jakubowski:86}, which considers
completely regular Hausdorff spaces (Tychonoff spaces)
of which $L_2(g_{\rho})$ with its weak topology is an example,
that a sufficient condition  is that 
\begin{enumerate}
\item
For each $T\ge 0$ and $\varepsilon >0$ there is a (weakly) compact subset
$K_{T,\varepsilon}$ of $L_2(g_{\rho})$ such that
\begin{equation}
\label{cpctcon}
\PP\left( Z^{N_j} \in \DD([0,T], K_{T,\varepsilon}) \right) > 1-\varepsilon\,,
\qquad j\ge1\,.
\end{equation}

\item
For each $d\ge1$, the 
$d$-dimensional processes $(Z^{N_j}(1), \ldots, Z^{N_j}(d))_{j\ge1}$
are tight.
\end{enumerate}

Lemma~\ref{tight.in} implies that the assumptions of Lemma~\ref{l2bds} hold,
and (\ref{cpctcon}) follows considering the Markov inequality.
We use (\ref{znavecg}) (derived from (\ref{zeq}))
 and (\ref{bracket}), and the  bounds (\ref{grz}), 
(\ref{decr1}) and (\ref{brkbd}). The uniform bounds in Lemma~\ref{l2bds}
and the fact that $Z^{N}(k)$ has jumps of size $N^{-1/2}$ imply
classically that $(Z^{N_j}(1), \ldots, Z^{N_j}(d))_{j\ge1}$
is tight, see for instance  Ethier-Kurtz~\cite{Ethier:86} Theorem~4.1 p.~354
or Joffe-M{\'e}tivier~\cite{Joffe:86} Proposition~3.2.3 
and their proofs.

\noindent{\em Step 2\/}.
The tightness result for $(Z^{N_j})_{j\ge1}$ 
implies it converges in law along some further subsequence to some $Z^\infty$
with initial law given by the law of $Z^\infty_0$. Considering (\ref{diff}), 
we have in (\ref{znavecg})
\begin{equation}
\label{fin}
N^{1/2}\left(F (R^N_s)(k) - F(\ut)(k) \right)
= \lint Z^N_s + N^{1/2} H\left(\ut,N^{-1/2}Z^N_s\right).
\end{equation}
We likewise consider (\ref{bracket}).
We use again the  bounds (\ref{grz}), (\ref{decr1}) and (\ref{brkbd}),
the uniform bounds in Lemma~\ref{l2bds},
and additionally (\ref{defH}) and Lemma~\ref{lin}. 
We deduce by a martingale characterization
that $Z^\infty$ has the law of the Ornstein-Uhlenbeck process
unique solution for (\ref{sde})
in $L_2(g_{\rho})$
starting at $Z^\infty_0$, see Theorem~\ref{Dbdd}. 
The drift vector is given by the limit for (\ref{zeq}) and
(\ref{znavecg}) considering (\ref{fin}),
and the diffusion matrix by the limit for (\ref{bracket}).
See for instance 
Ethier-Kurtz~\cite{Ethier:86} Theorem~4.1 p.~354 
or Joffe-M{\'e}tivier~\cite{Joffe:86} Theorem 3.3.1
and their proofs for details.

\noindent{\em Step 3\/}.
The limit in law of a sequence of stationary processes is stationary
(see Ethier-Kurtz~\cite{Ethier:86} p.~131, Lemma~7.7 and Theorem~7.8).
Hence the law of $Z^\infty$ is the unique law of the stationary 
Ornstein-Uhlenbeck process given by (\ref{sde}), see Theorem~\ref{ou.erg}.
 We deduce that from every subsequence we can extract a further
subsequence converging in law to this process.
Hence $(Z^N)_{N\ge L}$ converges in law
to this process.


\end{document}